\documentclass[12pt,oneside,a4paper]{amsart}
\usepackage{cmap}
\usepackage[T2A]{fontenc}
\usepackage[utf8]{inputenc}
\usepackage[ngerman, english]{babel}
\usepackage[pdftex,unicode]{hyperref}
\usepackage{amsmath}
\usepackage{amssymb}
\usepackage{amsthm}
\usepackage{verbatim}
\usepackage{relsize}
\usepackage{amsfonts}
\usepackage{graphicx}
\usepackage[normalem]{ulem}
\usepackage{extsizes}
\usepackage{float}
\usepackage{bbold}
\usepackage{dsfont}
\usepackage{calc}
\usepackage{bm}

\usepackage{fullpage}

\newlength{\widecommentlength}
\usepackage{pgfplots}
\newcommand{\widecommentbox}[3]{\def#1##1{\strut\newline\noindent\colorbox{#3}{\linespread{1}\parbox{.95\textwidth}{\small {\bf [#2]} ##1}}\newline}}

\renewcommand{\widecommentbox}[3]{\def#1##1{}}

\widecommentbox{\alex}{AP}{green!20!white}
\widecommentbox{\ad}{AD}{red!20!white}

\usepackage{enumitem}
\usepackage[square,numbers,sort&compress]{natbib}
\usepackage{mathtools}
\newcommand{\cal}[1]{\mathcal{#1}}
\renewcommand{\leq}{\leqslant}
\renewcommand{\geq}{\geqslant}
\renewcommand{\phi}{\varphi}
\newtheorem{theorem}{Theorem}
\newtheorem*{theorem*}{Theorem}
\newtheorem{prop}{Proposition}
\newtheorem{lemma}{Lemma}

\newtheorem*{corol*}{Corollary}
\theoremstyle{definition}
\newtheorem{definition}{Definition}
\newtheorem*{definition*}{Definition}
\newtheorem{example}{Example}
\theoremstyle{remark}
\newtheorem{remark}{Remark}
\newtheorem*{remark*}{Remark}

\newcommand\inner[2]{\langle #1, #2 \rangle}

\newcommand{\seq}[1]{\{{#1}_n\}_{n=1}^\infty}
\newcommand{\fsys}{\mathfrak{F}}
\newcommand{\fstarsys}{\mathfrak{F^{*}}}
\newcommand{\lattice}{\operatorname{Lat}\cal{A}}
\newcommand{\wt}{\hat{\cal{L}}}

\newcommand{\len}{\cal{L}}

\newcommand{\flow}{\mathcal{\hat{F}}}
\newcommand{\flowpos}{\mathcal{F}}
\newcommand{\preflow}{\mathcal{F^{*}}}
\newcommand{\flowposn}[1]{\mathcal{F}_{#1}}

\newcommand{\flowsgn}{\cal{\hat{F}}}
\newcommand{\source}{\mathbf{s}}
\newcommand{\sink}{\mathbf{t}}

\newcommand{\ter}{ter}
\newcommand{\ein}{in}
\newcommand{\eout}{out}
\newcommand{\eback}{\mathbf{back}}
\newcommand{\efor}{\mathbf{forward}}
\renewcommand{\root}{\mathbf{r}}

\newcommand{\net}{\Delta}
\newcommand{\onet}{\vec{\Delta}}
\newcommand{\gpaths}{\cal{P}_{G}}
\newcommand{\gnpaths}{\cal{P}_n}
\newcommand{\gstar}{G^{*}}
\newcommand{\gfi}{\varphi_{G}}

\newcommand{\phistar}{\phi_{*}}
\numberwithin{remark}{section}
\numberwithin{theorem}{section}
\numberwithin{prop}{section}
\numberwithin{equation}{section}
\numberwithin{lemma}{section}

\usepackage{caption}
\usepackage{tikz}
\usepackage{tkz-berge}
\usetikzlibrary{positioning, chains, fit, shapes, calc}
        \tikzset{lnode/.style={auto, circle, draw,fill=gray!50},
                 rnode/.style={auto, circle, draw},
                 mnode/.style={auto, circle, draw, fill=gray!25},
                 wnode/.style={sloped,above}}

\begin{document}
\title{Rank one density property for a class of~$M$-bases}
\author{Alexey Pyshkin}
\begin{abstract}
  In the early 1990s the works of Larson, Wogen and Argyros, Lambrou, Longstaff
    disclosed an example of a strong $M$-basis that did not admit a linear summation method.
  We study a class of $M$-bases $\mathfrak{F}=\{f_n\}_{n=1}^\infty$ in the Hilbert space that
    generalizes the Larson--Wogen system.
  We determine the conditions under which $\mathfrak{F}$ admits a linear summation method.
  In order to do that we employ some of the graph theory techniques.
\end{abstract}
\thanks{The author was supported by RFBR (the project~16-01-00674).}
\maketitle

\section{Introduction}
  Consider a complete minimal vector system $\fsys = \seq{f}$ in the infinite-dimensional Hilbert space $\cal{H}$.
  The sequence is called \emph{minimal} if none of its elements can be approximated by the linear combinations of the others.
  The system $\fsys$ is minimal when and only when it possesses a unique biorthogonal system $\fstarsys$.
  We say that $\fsys$ is an $M$-basis if $\fstarsys$ is complete as well.
  Let ${\cal{A} = \{T\in B(\cal{H}): Tf_n \in \overline{span}(f_n), \text{ for any } n \in \mathbb{N}\}}$ be an operator algebra, where
    $\overline{span}$ denotes a closed linear span;
    and $R_1(\cal{A})$ be an algebra generated by the rank one operators of $\cal{A}$.
  We denote by $\lattice$ the lattice of invariant subspaces for operators in $\cal{A}$.

  We are interested in the following three properties of the algebra $\cal{A}$.
  \begin{definition}[one point density property]
    \label{1pd}
    We say that $\fsys$ is \emph{one point dense} if for any $x \in \cal{H}$ and $\varepsilon > 0$
      there exists such $T\in \cal{A}$ that $||Tx - x|| < \varepsilon$.

  \end{definition}
    The definition is equivalent to $\fsys$ being a strong $M$-basis (see~\cite{katavolos}):
    the system $\fsys$ is called a strong $M$-basis if for any $x\in\cal{H}$ we have $x \in \overline{span}\big(\inner{x}{f^*_k}f_k\big)$, where
      $\overline{span}$ denotes a closed linear span.
  \begin{definition}[two point density property]
    \label{2pd}
    We say that $\fsys$ is \emph{two point dense} if for any $x, y \in \cal{H}$ and $\varepsilon > 0$
      there exists such $T\in \cal{A}$ that $||Tx - x|| < \varepsilon$ and $||Ty - y|| < \varepsilon$.
  \end{definition}
  \begin{definition}[rank one density property]
    \label{r1d}
    We say that the algebra $\cal{A}$ possesses \emph{rank one density property} if
      the unit ball of the rank one subalgebra $R_1(\cal{A})$ is dense in the unit ball of $\cal{A}$ in the strong operator topology.
  \end{definition}
  By abuse of notation, we say that the $\fsys$ is \emph{rank one dense}, meaning that the rank one density holds for the corresponding algebra $\cal{A}$.

  It is well known that the last definition is equivalent to $R_1(\cal{A})$ being dense
    in $\cal{A}$ in the ultraweak (or $\sigma$-weak) topology.
  Rank one density property can be considered as a generalization of the notion of \emph{linear summation method} for the system $\fsys$.
  We say that the system $\fsys$ admits a linear summation method if there exist $\alpha_{kn} \in \mathbb{C}, n \in \mathbb{N}, 0 < k < k_n$
    such that for any $x \in \cal{H}$,
    \[
      \lim_{n \to \infty} \sum_k \alpha_{kn} \inner{x}{f^*_k}f_k = x.
    \]
  Clearly, the existence of the linear summation method implies the rank one density property.

  Notice that the rank one density property~\eqref{r1d} implies the two point density property~\eqref{2pd}, and~\eqref{2pd}
    in its turn implies the one point density property~\eqref{1pd}.
  But is it true that $\eqref{1pd}$ implies $\eqref{r1d}$ or that~\eqref{2pd} implies~\eqref{r1d}?

  Actually, the first question has been studied for some time already.
  Long\-staff in~\cite{longstaff} studied the abstract subspace lattices and the corresponding operator algebras.
  Particularly, Longstaff proved that the rank one density property implies
    the complete distributivity of its lattice (see Raney~\cite{raney} for a definition).
  It is now a well-known fact that for a subspace lattice generated from an $M$-basis the complete distributivity of its lattice is equivalent to
    the property~\eqref{1pd} (see~\cite{argyroslambrou} for a proof).
  In the same paper the following question was raised: does the complete distributivity of $\lattice$ imply the rank one density property?
  The answer was already known to be positive in the case of a totally ordered lattice~\cite{erdos}, and
    Longstaff proved it for a finite-dimensional Hilbert space in~\cite{longstaff}.
  Later Laurie and Longstaff discovered that the rank one density holds in the case of commutative subspace lattices~\cite{laurielongstaff}.

  However, the solution for the general case remained unknown until Larson and Wogen showed that the answer is negative (see~\cite{larson}).
  They constructed an example of a vector system $\fsys$ such that the lattice $\lattice$ generated from $\fsys$
    was completely distributive, but the corresponding algebra $\cal{A}$
    did not possess the rank one density property.
  \begin{example}[Larson--Wogen system $\fsys_{LW}$]
    \label{lw-sys}
    For any $j > 0$ we define
    \begin{align*}
      &f_{2j-1}=-a_{2j-1}e_{2j-2} + e_{2j-1} + a_{2j}e_{2j}, \qquad &f_{2j}=e_{2j},\\
      &f^*_{2j}=-a_{2j}e_{2j-1}+e_{2j}+a_{2j+1}e_{2j+1}, \qquad &f^*_{2j-1}=e_{2j-1},
    \end{align*}
    where $a_k$ are nonzero complex numbers for any positive $k$ and $a_0 = 0$.
  \end{example}
  The construction presented by Larson and Wogen was remarkably simple and elementary, and was mentioned in several papers afterwards.
  For instance, this example was also studied in~\cite{argyroslambrou} (see Addendum) and by Azoff and Shehada in~\cite{azoff}, regarding
    the reflexivity of $\lattice$.
  Finally, Katavolos, Lambrou and Papadakis in~\cite{katavolos} performed a deep analysis of the density properties
    of this vector system and deduced that for this system the property~\eqref{1pd} does not imply~\eqref{r1d},
    however the system also admitted the property~\eqref{2pd}.
  We refer to this particular example as to Larson--Wogen vector system in this paper.

  It is natural to try to understand whether the two point density always implies the rank one density property for an arbitrary $M$-basis $f_k$.
  The aim of this paper is to study the rank one density property for a specific class of vector systems (which we call \emph{B-class})
    that represents a natural generalization of the Larson--Wogen vector system.

  In the next section we set up the notation and state the main theorem.

\bigskip
\section{B-class vector systems: rank one density problem}
    \label{fsys2graphs}
    We are interested in the rank one density property of the vector system $\fsys = \seq{f}$
      in a separable Hilbert space $\cal{H}$.
    Suppose that $\cal{H}$ has an orthonormal basis $\seq{e}$.
    Let $\fstarsys = \{f_n^{*}\}_{n=1}^\infty$ be a biorthogonal system to the original system $\seq{f}$.
    \begin{definition}
      We say that the vector system $\fsys$ belongs to the \emph{B-class} whenever the following conditions are satisfied:
      \begin{enumerate}[label=\textbf{C\arabic*}]
        \item \label{c1} both $\fsys$ and $\fstarsys$ are complete (or $\fsys$ is an $M$-basis);
        \item \label{c2} either $f_n = e_n$ or $f^*_n = e_n$ for any $n > 0$;
        \item \label{c3} $\inner{f_n}{e_n} = \inner{f^*_n}{e_n} = 1$ for any $n > 0$;
        \item \label{c4} $\inner{f_n}{e_k} = -\inner{f^*_k}{e_n}$ for any $n, k > 0, n \neq k$;
        \item \label{c5} the matrices $\{\inner{f_n}{e_k}\}$ and $\{\inner{f^*_n}{e_k}\}$ are both finite-band.
      \end{enumerate}
    \end{definition}
    \begin{prop}
      The definition given above guarantees the biorthogonality of the $f_n$ and $f^*_n$.
    \end{prop}
    \begin{proof}
      Consider the set of indices $N = \left\{n \in \mathbb{N} \mid f_n = e_n \right\}$.
      The property~\ref{c4} yields that $\inner{f_n}{e_k}$ is equal to zero whenever $k$ is not in
        $N \cup {n}$.
      Similarly, $f^*_m$ is orthogonal to $e_k$ for all $k$ in $N \setminus {m}$.
      We also have
      \[
        \inner{f_n}{f^*_m} = \sum_{k} \inner{f_n}{e_k} \inner{e_k}{f^*_m}.
      \]
      It follows from above that each summand on the right side is equal to zero whenever $k \neq n, m$.
      Thus we get
      \[
        \inner{f_n}{f^*_m} = \inner{f_n}{e_n} \inner{e_n}{f^*_m} + \inner{f_n}{e_m} \inner{e_m}{f^*_m}
        = \inner{e_n}{f^*_m} + \inner{f_n}{e_m} = 0,
      \]
      if $n \neq m$.
      Otherwise, if $n = m$, we get $\inner{f_n}{f^*_m} = 1$ due to the conditions~\ref{c2} and~\ref{c3}.
    \end{proof}
    The B-class of vector systems is the main subject of this section.
    It is a natural extension of the Larson--Wogen system (see Example~\ref{lw-sys}).

    Furthermore, a B-class vector system $\fsys$ could be associated with a
      locally-finite weighted bipartite graph $B(\fsys) = (V, E, \wt)$, where $\wt$ is a
      nonzero real-valued function on $V\times V$ such that $\wt(v, u) = -\wt(u, v)$.

    For each index $l > 0$ such that $f^*_l = e_l$ we put the vertex $v_l$ in the first part of the bipartite graph.
    We will call this part from now on the \emph{left} part of the graph $B(\fsys)$.
    For any other index $r > 0$ we construct a vertex in the other part of the graph.
    Evidently, for such indices $r > 0$ the condition $f_r = e_r$ holds due to the definition of the B-class vector systems.
    The second part of the graph will be referred as the \emph{right} part of the graph $B(\fsys)$.
    We put an edge between two vertices $v_l$ and $v_r$ from the left and right parts respectively,
      whenever the scalar product $\inner{f_l}{e_r}$ is not zero.
    For such two vertices we have $\inner{f_l}{e_r}$ = $-\inner{f_r}{e_l}$.
    We set the weight (length) on the edge $(lr)$,
    \[
      \wt_{lr} = \inner{f_l}{e_r}^{-1}.
    \]
    Obviously, $\wt_{lr} = -\wt_{rl}$.
    Since the Hilbert space is infinite-dimensional, the graph is infinite as well.
    However, due to the finite-band condition~\ref{c5}, the vertices of the graph are of a finite degree, hence the locally-finiteness of
      the constructed graph.
    \begin{remark}
      Note that there might be several vector systems for a single weighted bipartite graph $B$.
    \end{remark}
    \begin{remark}
      Observe that both parts of the constructed bipartite graph contain an infinite number of vertices, because
        the system $\fsys$ is an $M$-basis.
    \end{remark}
    You can see the bipartite graph built from the Larson--Wogen system in the Figure~\ref{lw-bgraph}.

    \bigskip
    Below we pose the main theorem.
    \begin{theorem}
      \label{thm-main}
      The B-class system $\fsys$ is rank one dense if and only if
        for any infinite path $\seq{r}$ in the bipartite graph $B(\fsys)$
          the series $\sum_{k=1}^\infty\lvert\wt(r_k, r_{k+1})\rvert$ diverges.
    \end{theorem}
    The bipartite graph $B(\fsys_{LW})$, associated with the Larson--Wogen system~\eqref{lw-sys},
      is a single ray itself.
    Therefore, we can apply Theorem~\ref{thm-main} to $\fsys_{LW}$, thus reestablishing
      the known fact from the papers~\cite{katavolos},~\cite{larson} and~\cite{argyroslambrou}.
    \begin{corol*}
      The system $\fsys_{LW}$ admits a rank one density property if and only if
        the sequence $\left\{1/a_n\right\}_{n=1}^\infty$ does not belong to $\ell^1$.
    \end{corol*}
    \begin{figure}
      \begin{center}
      \begin{tikzpicture}[thick,
                          every node/.style={draw},
                          ]

        \begin{scope}[start chain=going below,node distance=12mm]
          \foreach [evaluate={\k=int(\i*2 + 1)}] \i in {0,...,3}
          \node[lnode,on chain] (l\k) [label=left: $v_{\k}$] {$e_{\k}$};
          \node[draw=none,on chain, yshift=7mm]{$\vdots$};
        \end{scope}
        \begin{scope}[xshift=4cm,yshift=-0.7cm,start chain=going below,node distance=12mm]
          \foreach [evaluate={\k=int(\i*2 + 2)}] \i in {0,...,3}
          \node[rnode,on chain] (r\k) [label=right: $v_{\k}$] {$e_\k$};
          \node[draw=none,on chain, yshift=7mm]{$\vdots$};
        \end{scope}


        \begin{scope}[wnode]
          \foreach [evaluate={\l=int(\i*2 + 1);\r=int(\i*2 + 2)}] \i in {0,...,3}
          \draw (l\l) to node[draw=none, color=black] {\footnotesize $\mathrm{w}_{\l\r}=a_{\r}$} (r\r);
          \foreach [evaluate={\l=int(\i*2 + 3);\r=int(\i*2 + 2)}] \i in {0,...,2}
          \draw (l\l) to node[draw=none,color=black] {\footnotesize $\mathrm{w}_{\r\l}=a_{\l}$} (r\r);
        \end{scope}
      \end{tikzpicture}
      \caption{The bipartite graph $B(\fsys_{LW})$ } \label{lw-bgraph}
      \end{center}
    \end{figure}

    \medskip
    \subsection{Rank one density property for B-class}
      Suppose we have a B-class vector system $\fsys=\seq{f}$.
      First of all, we intend to demonstrate a reformulation of the rank one density problem for
        the B-class vector systems in terms of infinite networks.

      Recall that the system $\seq{f}$ is rank one dense
        if and only if there is no trace class operator $T: \cal{H} \to \cal{H}$ with the trace equal to $1$
        that belongs to the annihilator of the rank one subalgebra $R_1(\cal{A})$: one has $\inner{Tf_n}{f_n^*} = 0$ for any $n$
        (see Theorem 2.2 of~\cite{katavolos} for details).

      Suppose that there is an operator $T$ such that $\inner{Tf_n}{f_n^*} = 0$ for any $n$.
      There are two cases: either $f^*_n = e_n$ or $f_n = e_n$.
      In the first case the condition $\inner{Tf_n}{f_n^*} = 0$ turns into
      \begin{equation}
          \label{left-eqn}
          \sum_j T_{nj} \inner{f_n}{e_j} = 0,
      \end{equation}
      and in the second case it is equivalent to
      \begin{equation}
          \label{right-eqn}
          \sum_j T_{jn} \inner{f^*_n}{e_j} = 0.
      \end{equation}

      Now consider an auxiliary function $\flowsgn: V \times V \to \mathbb{R}$ defined as follows
      \begin{align*}
          &\flowsgn(v_l, v_r) = T_{lr} \inner{f_l}{e_r} = T_{lr} \wt_{lr}^{-1},\\
          &\flowsgn(v_r, v_l) = T_{lr} \inner{f^*_r}{e_l} = T_{lr} \wt_{rl}^{-1}.
      \end{align*}
      Observe that $\flowsgn$ is a skew-symmetric function.
      Moreover, two equalities~\eqref{left-eqn} and~\eqref{right-eqn} correspond to the \emph{left} and \emph{right}
        parts of the bipartite graph $B(\fsys)$ respectively.
      It follows that the condition $\inner{Tf_n}{f_n^*} = 0$ could be reduced to a simpler one:
      \begin{equation}
        \label{almost-flow-eqn}
        \sum_{u \in V} \flowsgn(v, u) + T_{vv} = 0
      \end{equation}
        for each vertex $v$ in the graph $B(\fsys)$.
      \begin{remark}
        Observe that the function $\flowsgn$ defined on the graph $B(\fsys)$ resembles
          a \emph{flow} defined on the edges of the graph $B(\fsys)$.
        One might also see that the equation~\eqref{almost-flow-eqn} describes the total flow (sum of the outgoing flows and incoming flows) 
          for each vertex $v$ in the graph $B(\fsys)$.
        In order to formalize this observation we are going to build a flow
          after a few changes are made to the graph $B(\fsys)$.
      \end{remark}

    \medskip
    \subsection{Flows and networks preliminaries}

      Before we proceed we are going to introduce a few basic definitions.
      \begin{definition}
        \emph{Network} $\net$ is a quadruple $(G, \len, \source, \sink)$, where $G = (V, E)$ is a weighted graph
        with a positive length function $\len$ on $E$ and two vertices $\source, \sink \in V$, which
        we will call \emph{source} and \emph{sink} of the network respectively.
      \end{definition}
      By graph here we always mean a graph without loops and multiple edges.
      \begin{definition}
        The skew-symmetric function $\flow: V \times V \to \mathbb{R}$ is called a \emph{pseudo-flow}.
      \end{definition}
      \begin{definition}
        Let $G = (V, \vec{E})$ be an oriented graph.
        For each vertex $v \in V$ we take $\ein_G(v)$ as the set of incoming edges
          and $\eout_G(v)$ as the set of outgoing edges in the graph $G$.
        We will omit the graph from the notation whenever it is clear from the context.
        For the set of vertices $V_0$ we denote
        \begin{align*}
          &\ein(V_0) = \big\{(uv)\in E \mid u \in V \setminus V_0, v \in V_0\big\},\\
          &\eout(V_0) = \big\{(vu) \in E \mid u \in V\setminus V_0, v \in V_0\big\}.
        \end{align*}
        For the subgraph $G_0 = (V_0, \vec{E}_0) \subseteq G$ we use the same notation.
        \begin{align*}
          &\ein(G_0) = \big\{(uv)\in \vec{E}_0 \mid u \in V \setminus V_0, v \in V_0\big\},\\
          &\eout(G_0) = \big\{(vu) \in \vec{E}_0 \mid u \in V\setminus V_0, v \in V_0\big\}.
        \end{align*}
        For each edge $e=(uv) \in \vec{E}$ we will write 
          $\ter(e)$ for the ending vertex $v$.
        Set $\ter(\vec{E}_0) = \big\{\ter(e) \mid e\in \vec{E}_0\big\}$.
      \end{definition}
      \begin{definition}
        Let $\net = (G, \len, \source, \sink)$ be a network, and $\flow$ be a pseudo-flow.
        Then $d_{+}(v)$ will stand for the sum of the pseudo-flows \emph{leaving} the vertex $v$ and
        $d^{-}(v)$ will stand for the sum of the pseudo-flows \emph{entering} $v$:
        \begin{align*}
          &d^{+}_{\flow}(v) = \sum_{\flow(vu) > 0} \flow(v,u),\\
          &d^{-}_{\flow}(v) = \sum_{\flow(uv) > 0} \flow(u,v).
        \end{align*}
        Let $d_{\flow}(v)$ be equal to $d^{+}_{\flow}(v) - d^{-}_{\flow}(v)$.
        We will refer to this value as a \emph{total flow} of the vertex $v$.
        The vertex $v$ is called \emph{$\flow$-active} if $d(v)$ is less than zero,
          \emph{$\flow$-deficient} if $d(v)$ is greater than zero and
          \emph{$\flow$-preserving} if $d(v)$ is precisely zero, meaning that the total incoming flow
          is equal to the total outgoing flow of the vertex $v$.
      \end{definition}
      \begin{definition}
        \label{flow-dfn}
        Given a network $\net = (V, E, \len, \source, \sink)$ and a pseudo-flow $\flow$ we will call $\flow$ a
          \emph{flow} if for any $v \in V$ the total flow is correctly defined (meaning that the corresponding sum converges absolutely),
          and for any $v \in V \setminus \{\source, \sink\}$ the total flow is zero: $d(v) = 0$.
      \end{definition}
      \begin{definition}
        \label{flow-preserving-dfn}
        For a network $\net = (G, \len, \source, \sink)$ and a flow $\flow$ we will say
          that the network $\net$ \emph{preserves} the flow $\flow$ if the total flows $d(\sink)$, $d(\source)$ are finite
          and $d(\sink) = -d(\source)$.
      \end{definition}
      \begin{remark}
        In simple words this property suggests that the total flow coming out of the source is equal to the total flow
          coming into the sink.
        Note that in the case when the graph $G$ is finite, the network $\net$ always preserves the flow.
        It is the infinite case that is of interest.
      \end{remark}
      \begin{definition}
        \emph{Oriented network} $\onet$ is a quadruple $(\vec{G}, \len, \source, \sink)$, where
          $\vec{G} = (V, \vec{E})$ is an oriented weighted graph with
          a positive length function $\len$ on $\vec{E}$ and two vertices $\source, \sink \in V$, which
          we will call the \emph{source} and the \emph{sink} of the oriented network respectively.
      \end{definition}
      Now we define the flow functions on the oriented networks in the similar manner we defined on the non-oriented networks.
      \begin{align*}
        &d^{+}_{\flowpos}(v) = \sum_{e \in \eout(v)} \flowpos(e),\\
        &d^{-}_{\flowpos}(v) = \sum_{e \in \ein(v)} \flowpos(e),\\
        &d_{\flowpos}(v) = d^{+}_\flowpos(v) - d^{-}_\flowpos(v).
      \end{align*}
      Sometimes we will omit the flow from the notation.
      The definitions~\ref{flow-dfn} and~\ref{flow-preserving-dfn} are the same for the oriented networks.
      \bigskip

      Suppose $\net = (V, E, \len, \source, \sink)$ is a network and $\flow$ is a skew-symmetric function defined
        on $\net$.
      Naturally, we might instead consider an oriented network $\onet = (V, \vec{E}, \len, \source, \sink)$ and a
        positive flow $\flowpos: \vec{E} \to \mathbb{R}^{+}$, where $V$, $\source$, $\sink$ and $\len$ are the same.
      We assign a direction to the edge $(uv)$ as follows:
        $(uv) \in \vec{E}$ if $\flow(u, v) > 0$, and $(vu) \in \vec{E}$ otherwise.
      This way only one of the edges $(uv)$ and $(vu)$ is present in the graph $\vec{G} = (V, \vec{E})$.
      Also we set the positive flow $\flowpos: E \to \mathbb{R}^{+}$ so that
        for any edge $(uv) \in \vec{E}$ we have $\flowpos(u v) = \flow(u,v)$.
      \begin{remark}
        The flow functions on the networks and the positive flow functions on the oriented networks
          are interchangeable and describe the same object.
        The notion of the direction of the positive flow over the particular edge $(uv)$ of the oriented network
          is incorporated into the sign of the flow on the same edge $(uv)$ of the non-oriented network.
      \end{remark}

    \medskip
    \subsection{B-network construction}
      \begin{definition}
        Consider a network $\net = (G, \len, \source, \sink)$ such that
          the degree of each vertex in $V \setminus \{\source, \sink\} $ is finite,
          and the length of each edge incident to source or to sink is equal to one.
        Also we demand that the vertices $\source$ and $\sink$ are connected by a finite path.
        We will call such network a \emph{B-network}.
      \end{definition}

      In this section we build up a B-network from the graph $B(\fsys)$ and associate a flow $\flow$ on this network with 
        an arbitrary trace class operator $T$ which annihilates the rank one subalgebra $R_1(\cal{A})$.
      We will see later that this association is a bijection.
      Namely, we plan to construct the network $\net(\fsys) = (V, E, \len, \source, \sink)$,
        with the length function $\len: E \to \mathbb{R}^{+}$ and a real skew-symmetric flow $\flow: V \times V \to \mathbb{R}$.
      Firstly, we incorporate all the vertices and edges from the graph $B(\fsys)$ into the network $\net(\fsys)$.
      We set the flow and the length functions on the edge $e = (v_l v_r)$ as follows:
      \begin{align*}
        &\flow(v_l, v_r) = T_{lr} \inner{f^*_r}{e_l},\\
        &\flow(v_r, v_l) = T_{lr} \inner{f_l}{e_r} = -\flow(v_l, v_r),\\
        &\len(e) = \lvert \wt(v_l, v_r) \rvert =  \lvert\inner{f_l}{e_r}\rvert^{-1}.
      \end{align*}
      \begin{remark}
        Evidently on this kind of edges the flow $\flow$ agrees with the function $\flowsgn$
          we examined a few paragraphs before.
      \end{remark}

      Now we add two new vertices: the source vertex $\source$ and the sink vertex $\sink$ to the constructed graph.
      For each vertex $v_l$ from the left part of the graph we connect it with the vertex $\source$ with the edge $e_l = (\source v_l)$
        and assign the flow to the newly constructed edge.
      \begin{align*}
        &\flow(v_l, \source) = -\flow(\source, v_l) = T_{ll}.
      \end{align*}
      We set the length equal to one for such edges: $\len(e_l) = 1$.
      As we added the edges $\seq{e}$, the flow became \emph{preserved} at each vertex of the left part of the network $\net(\fsys)$.

      Likewise, for any vertex $v_r$ from the right part of the graph we add a new edge $e_r=(v_r \sink)$
        and set the flow $\flow$ equal to:
      \begin{align*}
        &\flow(v_r, \sink) = -\flow(\sink, v_r) =  T_{rr}.
      \end{align*}
      Again, we have $\len(e_r) = 1$ for any $k$.

      The network $\net(\fsys)$ is obviously a B-network.
      The defined function $\flow$ is a flow, since for vertex the total flow is correctly defined,
        and each of the vertices preserves $\flow$ due to~\eqref{almost-flow-eqn}.
      \begin{remark}
        Note that the network $\net(\fsys)$ depends only on the biorthogonal system $\fsys$ and not on the operator $T$.
        Only the flow $\flow$ depends on the operator $T$.
      \end{remark}
      What can be said about the total flow in each of the vertex of the constructed network $\net(\fsys)$?
      Due to the trick we performed, the total flow became zero in each of the vertices from the left and the right parts.
      The total flow $d(\source)$ in the source vertex is equal to $\sum T_{ll}$, and the total flow in the
        sink vertex is now equal to $\sum T_{rr}$.
      Now one can see that the network $\net(\fsys)$ is $\flow$-preserving if and only if the trace of the operator $T$ is equal to zero.
      \begin{definition}
        Consider a network $\net = (V, E, \len, \source, \sink)$ and a flow $\flow$ on it.
        We write $\lvert\flow\rvert$ for the \emph{mass} of the flow $\flow$:
        \[
          \lvert\flow\rvert = \sum_{e \in E} \lvert\flow(e) \len(e)\rvert.
        \]
        For the flow $\flowpos$ on the oriented network we will use the same notation.
      \end{definition}
      \begin{theorem}
        \label{thm-graph-eq}
          Let $\net(\fsys) = B(\fsys, \len, \source, \sink)$ be a B-network constructed
            from the B-class system $\fsys = \seq{f}$.
          Then $\fsys$ is rank one dense if and only if
            the network $\net(\fsys)$ is $\flow$-preserving for any finite-mass flow $\flow$.
      \end{theorem}
      \begin{proof}
        Suppose there is no rank one density property for the system $\seq{f}$.
        As we mentioned before, it implies that there exists a trace class operator $T: \cal{H} \to \cal{H}$ with
          the trace equal to one, such that $\inner{Tf_n}{f_n^*} = 0$ for any $n$.
        Using the operator $T$, we were able to define a flow $\flow$ on the network $\net(\fsys)$.
        The constructed flow has a finite mass because the operator matrix of $T$ is finite-band and $T$ has a finite trace.
        Finally, the network $\net(\fsys)$ does not preserve the flow since
        \[
          d(\source) + d(\sink) = \sum_{l} T_{ll} + \sum_{r} T_{rr} = Tr(T) = 1,
        \]
        thus $d(\sink) \neq -d(\source)$.
        The necessity is proved.

        Suppose that the system $\fsys$ is rank one dense.
        Assume that there is also a flow $\flow$ on the network $\net(\fsys)$ such that $d(\sink) + d(\source) \neq 0$.
        Firstly, consider the edges incident to the source vertex.
        Having all the vertices already enumerated we will consider the vertices in the left part $v_{l}$.
        Recall that each vertex $v_i$ (except source and sink) matches to the basis element $e_i$.
        Assign the diagonal elements of $T$:
        \[
          T_{ll} = \flow(v_l, \source) \len(v_l v_l) = \flow(v_l, \source), \quad l > 0.
        \]
        We proceed in this fashion with the right part of the graph $\net(\fsys)$:
        \[
          T_{rr} = \flow(v_r, \sink) \len(v_r v_r) = \flow(v_r, \sink), \quad r > 0.
        \]
        Now consider two connected vertices $v_l$, $v_r$ from the left and right part respectively.
        Let $T_{lr} = \flow(v_l, v_r) \len(v_l  v_r)$.
        Let all the other matrix elements of the operator matrix $T_{ij}$ be zero.
        Then observe that we obtained a finite-band operator matrix $T_{ij}$ with the sum of diagonal elements equal to one.
        Since the mass of the flow $\flow$ is finite, $T_{ij}$ is a summable sequence, considering that
          the $\ell^1$-norm of $T_{ij}$ is exactly the mass of $\flow$.
        Then the operator $T$ has a finite nonzero trace, since $Tr(T) = d(\sink) + d(\source)$, which is not zero
          by the assumption we made.
        Due to the fact~\eqref{almost-flow-eqn} that $\flow$ is preserved at each vertex of the network $\net(\fsys)$,
          we can see that the operator $T$ annihilates all the rank one operators $R_1(\cal{A})$,
          which is a contradiction to the rank one density property.
      \end{proof}

  \bigskip
  \section{B-networks characterisation}
    Due to Theorem~\ref{thm-graph-eq}, we are able to analyze the flows on the B-network
      $\net(\fsys)$ in order to understand the conditions under which $\fsys$ is rank one dense.
    In this section we are not going to address the Hilbert space setup at all,
      instead we will use only the abstract objects of the graph theory we introduced in the previous section.
    It will be more convenient for us to study the \emph{oriented} networks in this section unlike we
      did in the previous one.
    Consequently, the graphs, edges and networks are presumed to be oriented unless said otherwise.

    Next we introduce a few more notions from the graph theory.
    \begin{definition}
      Let $G = (V, E)$ be an oriented graph.
      A sequence of vertices $\{v_k\}_{k=1}^N$ is called a \emph{path} if each pair of consecutive vertices
        $v_k$ and $v_{k+1}$ are connected by an edge $(v_k v_{k+1})$ and no edges or vertices are repeated twice in the sequence.
      The path is called a \emph{ray} when $N$ is equal to infinity.
    \end{definition}
    \begin{definition}
      Let $G$ be an oriented graph.
      A finite sequence of vertices $\{v_k\}_{k=1}^N$ is called a \emph{cycle} (or a simple cycle) if it is
        a path and there is an edge $e = (v_N v_1)$ connecting the last vertex to the first vertex of the sequence.
    \end{definition}

    The main result of this section is Theorem~\ref{ray-thm}.
    \begin{theorem}
      \label{ray-thm}
      Let $\net = (G, \len, \source, \sink)$ be an oriented B-network.
      Then a positive flow $\flowpos: E \to \mathbb{R}^{+}$ such that
      \begin{enumerate}[label=(\alph*)]
        \item the flow $\flowpos$ has a finite mass,
        \item the network $\net$ does not preserve the flow,
      \end{enumerate}
        exists if and only if there exists a ray $r = \seq{r}$ in the network $\net$, taking its start in the vertex $\source$,
          with a finite total length:
        \[
          \sum_{k=1}^\infty \len(r_k r_{k+1}) < \infty.
        \]
    \end{theorem}
    \begin{proof}
      Without loss of generality we might assume that $d(\sink) + d(\source) = 1$.
      To simplify the proof, we consider a modified network $\net' = (G', \source', \source')$.
      Namely, we merge the source vertex $\source$ with the sink vertex $\sink$ and call this a new source vertex $\source'$.
      This way we obtain a network $\net'$ with a single source such that $d(\source') = 1$.
      In that case the graph $G'$ is a rooted graph with the root at the vertex $\source'$.
      From now on we will employ the term \emph{root vertex} and denote it by $\root = \source$.

      The flow $\flowpos$ as well as the length function $\len$ stay exactly the same (being defined on the edges set $E'$).
      All the vertices in the graph continue preserving the flow after this procedure, except for the source vertex $\source'$,
        which is now \emph{$\flowpos$-active} since its total flow $d(\source')$ is equal to one.
      Thus $\net'$ is a network with a function $\flowpos$ being a flow as required.
      Obviously, the mass of the flow does not change since we did not change the lengths of the edges.
      Furthermore, the network $\net'$ does not preserve the flow $\flowpos$ by construction.

      It is clear that $\net'$ contains a finite length ray starting on the vertex $s'$ if and only if the original network $\net$ contains a finite
        length ray starting on the vertex $s$.

      Therefore, we will be proving the following version of Theorem~\ref{ray-thm}.
      \begin{theorem*}
        Let $\net = (G, \len, \source, \sink)$ be an oriented B-network such that $\source = \sink$.
        Then a positive finite-mass flow $\flowpos: E \to \mathbb{R}^{+}$ such that
          $\net$ does not preserve $\flowpos$ (it means that $d(\source) \neq 0$)
          exists if and only if there exists a ray $r = \seq{r}$ which has a finite total length and starts on $\source$.
      \end{theorem*}
      First we prove the necessity.
      Suppose that there is a ray of a finite total length $r=(r_1 r_2 \dots)$.
      In that case we set the resulting flow $\flowpos(r_k r_{k+1})$ to one, creating a flow which escapes
        from $\source$ to infinity.
      All the other edges will not carry any flow, hence the flow is preserved
        at each vertex but not preserved by the whole network $\net$.
      The necessity is proved, now we turn to the sufficiency.

      \medskip
      \subsection{Eliminating the positive flow cycles}
        For a start we get rid of all the positive flow cycles in the oriented graph $G$.
        \emph{A positive flow cycle} is a cycle $C$ in $G$ such that positive flow circulates along $C$,
          meaning that for each edge $e \in C$ we have $\flowpos(e) > 0$.
        In order to accomplish that, we will present a sequence of the flows $\flowposn{n}$, defined on the network $\net$.
        At each step $n \in \mathbb{N}$ we examine a subgraph $G_n \subset G$ and present a flow $\flowposn{n}: \net \to \mathbb{R}^{+}$
          such that $\flowposn{n}$ is monotone: for any edge $e \in E$ we have $\flowposn{n+1}(e) \leq \flowposn{n}(e)$.
        The main property of the flow $\flowposn{n}$ is that there will be no positive flow cycles in the induced subgraph
        $G_n$.
        We start with the flow $\flowposn{0}:\net \to \mathbb{R}^+$ equal to the flow $\flowpos$.
        The first step $n = 1$ is trivial: the graph $G_1$ contains only the root vertex, and
          the flow $\flowposn{1}$ is equal to $\flowposn{0}$ at any edge $e \in E$.

        In order to proceed we enumerate all vertices in the network $V = \left\{v_k\right\}_{k=1}^\infty$.
        At the step $n$ we consider the subgraph $G_n = (V_n, E_n) \subset G$, where $V_n$ is defined as $\left\{v_k\right\}_{k=1}^n$
          and the edge set $E_n$ consists of the edges from $E$ incident to $V_n$.

        Suppose we have already established the flow $\flowposn{n}: \net \to \mathbb{R}^{+}$ for which
          all the inequalities above are satisfied.
        We are aiming to construct the flow $\flowposn{n+1}$.

        Define the flow function $\cal{E}_0(e): \net \to \mathbb{R}^{+}$ as $\cal{E}_0(e) = \flowposn{n}(e).$
        Let there be a simple cycle $C_1$ in the graph $G_n$ such that each edge $e \in C_1$ carries
          a positive flow: $\cal{E}_0(e) > 0$.
        At this point we can decrease the flow $\cal{E}_0$ on the edges of $C_1$ so that the cycle $C_1$
          ceases being a positive flow cycle, total mass of the flow does not increase, and all the other properties stay intact.
        Denote by $\min(C_1)$ the minimal value the flow $\cal{E}_0$ attains on the edges of $C_1$.
        Consider an adjusted flow $\cal{E}_1$.
        \begin{equation*}
          \cal{E}_1(e) =
          \begin{cases}
            \cal{E}_0(e) - \min(C_1) \quad \text{if } e \in C_1,\\
            \cal{E}_0(e) \quad \text{otherwise}.
          \end{cases}
        \end{equation*}
        Since $C_1$ is a cycle, we reduced the flow $\cal{E}_0$ in each edge by the same value,
          so $\cal{E}_1$ is indeed a flow.
        After this procedure the positive flow cycle $C_1$ disappears,
          leaving us with a lesser flow than we had before this step: $\cal{E}_1(e) \leq \cal{E}_0(e)$ for any $e \in E_n.$
        In this manner we remove all the positive flow cycles from the finite graph $G_n$.
        It is achievable since at each step we turn at least one edge into a zero-flow edge (such edge $e$ that $\cal{E}_0(e) = 0$).

        Suppose we removed $K$ positive flow cycles in total from the graph $G_n$.
        Then set the flow $\flowposn{n+1}$ equal to the final flow $\cal{E}_K$: $\flowposn{n+1} = \cal{E}_K.$
        We have $\flowposn{n+1}(e) \leq \flowposn{n}(e) \leq \flowpos(e)$ for each edge $e \in E$.

        Now set $\preflow(e) = \lim_{n \to \infty} \flowposn{n}(e)$ for each $e \in E$.
        The described pointwise limit exists since for each edge $e$ the sequence $\flowposn{n}(e)$ is decreasing and bounded below.
        \begin{prop}
          The preflow $\preflow$ is a flow such that $d_{\preflow}(\root) = d_{\flowpos}(\root)$.
          The graph $G$ contains no positive flow cycles of $\preflow$.
        \end{prop}
        \begin{proof}
          Obviously, for each $n$ the function $\flowposn{n}$ is a flow since we do not affect the total flow each time we remove a positive flow cycle.

          Now we need to obtain the same for the limit case.
          Firstly, look at the non-root vertex $v \neq \root$.
          Since there are a finite number of edges in $E$ incident with $v$, the limit of
            a finite sum $d_{\flowposn{n}}(v)$ is equal to $d_{\preflow}(v)$.
          The vertex $v$ is not a root, so $d_{\flowposn{n}}(v) = 0$ for any $n$, and so the limit $d_{\preflow}(v)$ is equal to zero as well.
          Consequently, the function $\preflow$ is a flow in the network $\net$.

          Secondly, consider the root vertex $\root$.
          Notice that
          \[
            d_{\flowposn{n}}(\root) \overset{\mathit{def}}{=} \sum_{\eout(\root)} \flowposn{n}(e) - \sum_{\ein(\root)} \flowposn{n}(e).
          \]
          Since the sequence $\{\flowposn{n}(e)\}$ is decreasing and bounded from below,
            the quantity $d_{\preflow}(\root)=\sum_{\eout(\root)} \preflow(e) - \sum_{\ein(\root)} \preflow(e)$ is defined correctly.
          Because $d_{\flowposn{n}}(\root)$ is equal to $d_\flowpos(\root)$ for any $n$ due to the given algorithm, we have
          \[
            d_\preflow(\root) = \lim_{n \to \infty} d_{\flowposn{n}}(\root) = d_\flowpos(\root).
          \]
          Thus, $\preflow$ is a flow, let us prove the second part of the proposition.
          To obtain a contradiction, suppose that there is a cycle $C$ such that $\preflow(e) > 0$ for
            each edge $e \in C$.
          Clearly, there exists such $n$ that the cycle $C$ lies within the graph $G_n$.
          Since at each step we reduced the flow $\flowposn{k}$, at the step $n$
            the cycle $C$ was a positive flow cycle as well: $\flowposn{n}(e) > 0$ for each edge $e \in C$.
          That is a contradiction.
        \end{proof}

        Suppose that the flow $\flowpos$ contains no positive flow cycles in the first place.
        Also we remove the edges of the network $\net$ for which the flow $\flowpos$ equals zero.

      \medskip
      \subsection{Subgraphs construction}
        We are going to need a few definitions in order to proceed with the proof of the theorem.
        \begin{definition}
          Let $G = (V, E,\len, \root)$ be a weighted directed graph with a positive
            length function $\len$ defined on $E$ and a root vertex $\root \in V$.
          Then $\gpaths(u, v)$ denotes the set of paths from the vertex $u$ to the vertex $v$ in the graph $G$.
          Sometimes we will omit the first argument, and in that case we will refer to the paths from the root vertex
            $\gpaths(v) = \gpaths(\root, v)$.
          Furthermore, we expand this notion onto the vertices sets: $\gpaths(U,V)$ is the set of paths
            from the set of vertices $U$ to the set of vertices $V$ in the graph $G$.
          We write $\gpaths$ for the set of all paths in the graph $G$.
        \end{definition}
        \begin{definition}
          The distance function $\gfi: \gpaths \to \mathbb{R}^{+}$ is given by $\gfi(p) = \sum_{e \in p} \len(e)$.
          We extend this function to take values on the set of the vertices $V$ as well,
            $\gfi(v) = \displaystyle\inf_{\gpaths(v)} \gfi(p).$
        \end{definition}
        \begin{remark}
          Our ambition in the following paragraphs is to provide an estimate on the values of the distance function $\gfi$ on any vertex $v$
            located ''far'' from $\root$, namely for any integer $N > 0$ we are to provide $M > 0$ such that
            for any $v$ satisfying $\displaystyle\inf_{p \in \gpaths(v)} \#\left\{ e \mid e \in p\right\} > N$, one has $\gfi(v) < M$.
        \end{remark}

        The construction is based on the breadth-first search in the infinite graph $G$.
        We are going to define \emph{finite} subgraphs $G_n = (V_n, E_n)$ of the network $\net$ such that $G_{n-1} \subset G_{n}$ for each $n > 1$.
        We denote by $L_n$ the vertex set difference $V_n \setminus V_{n-1}$.
        Besides, we are going to establish a sequence of the positive preflow functions $\flowposn{n}: E_n \to \mathbb{R}^{+}$,
          pertaining the following properties.
        \begin{enumerate}[label=\textbf{P\arabic*}]
          \item \label{p1} $\flowposn{n+1}(e) \geq \flowposn{n}(e)$ for each $e \in E_n$;
          \item \label{p2} one has $\flowposn{n}(e) \leq \flowpos(e)$ for each $e \in E_n$;
          \item \label{p3} for the root vertex one has $d_{\flowposn{n}}(\root) = 1$;
          \item \label{p4} for any non-root $v \in G_{n-1}$ one has $d_{\flowposn{n}}(v) = 0$, for any $n > 1$;
          \item \label{p5} for each vertex $v \in L_n$ one has $d_{\flowposn{n}}(v) = -d^{-}_{\flowposn{n}}(v) < 0$.
        \end{enumerate}
        First of all, we scale the flow function $\flowpos$ at the root vertex so that $d^{+}_{\flowpos}(\root) > 1$.
        We are going to need this assumption at the first step.
      \medskip
      \subsection{Constructing $G_n$ and $L_n$}
        Let $G_0$ be the trivial subgraph containing the root vertex only: $G_0 = \left\{\root\right\}$.
        At the first step we build the graph $G_1$.
        We choose the set of the edges $E_1 \subseteq \eout(\root)$ such that $\sum_{e \in E_1} \flowpos(e) \geq 1$.
        It is possible since $d^{+}_{\flowpos}(\root) > 1$.

        Then we define the \emph{first layer} $L_1$ as the end vertices of $E_1$: $L_1 = \ter(E_1)$.
        Along with the graph $G_1 = (V_0 \cup L_1, E_1)$ we set the preflow $\flowposn{1}: E_1 \to \mathbb{R}^{+}$
        \[
          \flowposn{1}(e) = \frac{\flowpos(e)}{\sum_{e\in E_1}{\flowpos(e)}}, \quad e \in E_1.
        \]
        \begin{prop}
          The preflow $\flowposn{1}$ satisfies all of the properties~\ref{p1}--\ref{p5}.
        \end{prop}
        \begin{proof}
          The properties~\ref{p1}--\ref{p4} are trivial to check.

          We included the ends of the $E_1$ edges without any additional edges.
          Thus each vertex in $L_1$ has only incoming edges, which makes their
            total flow strictly negative (recall that we removed all zero-flow edges).
          Thus~\ref{p5} holds as well.
        \end{proof}
        \begin{remark}
          We dealt with the possibly infinite degree of the root vertex by abandoning some of the children of the root.
          \alex{a figure of the first layer (infinite degree!)}
          We need $G_n$ to be finite in order to analyze the function $\gfi$.
        \end{remark}

        We proceed by induction.
        Suppose that $G_{n-1}$, $L_{n-1}$ and the preflow $\flowposn{n-1}$ are already defined.
        Consider the outgoing edges of the vertex $v \in L_{n-1}$.
        Some of them might lead to some new vertices, which do not belong to $G_{n-1}$.
        Such edges will be referred to as $\efor_n$,
        \[
          \efor_n = \Big\{ e = (uv) \mid u \in L_{n-1} \text{ and } v \notin G_{n-1}\Big\}.
        \]
        Other edges might lead to the already visited vertices of $G$, namely the ones residing in the graph $G_{n-1}$.
        Such edges we will address as the \emph{back} edges of $L_{n-1}$ and denote them by $\eback_n$,
        \[
          \eback_n = \Big\{e = (us) \mid u \in L_{n-1} \text{ and } s \in G_{n-1}\Big\}.
        \]
        \begin{remark}
          Note that \emph{back} edge could never lead to the root vertex since that
            would yield an existence of a positive flow cycle for the flow $\flowpos$.
        \end{remark}
        \alex{a figure of the $n$-th layer construction}
        Finally, we set $G_n$ and $L_n$.
        \begin{align*}
          &L_n = \ter(\efor_n),\\
          &G_n = (V_{n-1} \cup L_n, E_{n-1} \cup \eback_n \cup \efor_n).
        \end{align*}

      \medskip
      \subsection{$n$-th step: setting a preflow}
        We divide the procedure in two steps.
        Consider the preflow function $\cal{G}: E_n \to \mathbb{R}^{+}$.
        \begin{align*}
          &\cal{G}(e) = \flowposn{n-1}, \quad e \in E_{n-1},\\
          &\cal{G}(e) = 0, \quad e \in E_n \setminus E_{n-1}.
        \end{align*}
        We are going to extend it further onto $E_n$ (step one) and then change it gradually in order
          to satisfy all the properties~\ref{p1}--~\ref{p5} (step two).
        Notice that all the vertices except the root and the last layer $L_{n-1}$ are $\cal{G}$-preserving.
        The root is $\cal{G}$-deficient and the $(n-1)$-th layer is currently $\cal{G}$-active.
        We first change the flow in the vertices $L_{n-1}$ in order to make them preserve the flow $\cal{G}$.
        For this purpose, we propagate the flow $\cal{G}$, incoming to the $(n-1)$-th layer, one edge further.

        Namely, consider the vertex $u \in L_{n-1}$.
        By~\ref{p2} we have $d_{\cal{G}}^{-}(u) \leq d_{\flowpos}^{-}(u)$,
          so we are able to set the flow $\cal{G}$ on the outgoing edges of $u$ in the following way:
        \begin{enumerate}[label=\textbf{(\roman*)}]
          \item\label{posprop} each outgoing edge has a positive flow: $\cal{G}(e) > 0$ for any $e \in \eout(u)$;
          \item $\cal{G}$ is still bounded from above: $\cal{G}(e) \leq \flowpos(e)$ for any $e \in \eout(u)$;
          \item the vertex $u \in L_{n-1}$ becomes $\cal{G}$-preserving.
        \end{enumerate}
        \begin{remark}
          Literally, here we are pushing the flow out of the vertex $u \in L_{n-1}$ and spreading it among all the outgoing edges
            $\eout_{G_n}(u) = \eout_G(u)$.
          Since the preflow $\cal{G}$ is bounded above by the flow $\flowpos$, which is preserved at each vertex of the graph $G$,
            particularly, at the vertex $u$, it is trivial that such an expansion exists.
        \end{remark}
        So, if we had no back edges, the function $\cal{G}$ would be the next preflow function $\flowposn{n}$, and this step would be accomplished,
          since all the properties~\ref{p1}--\ref{p5} would be complied by $\cal{G}$.
        However, the flow pushed from the $(n-1)$-th layer across the back edges made some of the vertices in $G_{n-1}$ \emph{not} flow-preserving.
        Namely, the end vertices of such edges $\eback_{n}$ now became $\cal{G}$-active, since for such vertices
          the total incoming flow now exceeds the total outgoing flow.
        \alex{a figure of a back edge making some inner vertex non-preserving, unite with the figure of n-th layer?}
        Finally, we have reached the step number two of the $\flowposn{n}$ construction.
        We are going to carry out the \emph{relaxation} procedure in order to alleviate the flow excess in those $\cal{G}$-active vertices.

        For simplicity of notation we continue to write $\gnpaths$ instead of $\cal{P}_{G_n}$.
        Consider a $\cal{G}$-active vertex $s \in \ter(\eback_n)$ such that $d_{\cal{G}}^{-}(s) > d_{\cal{G}}^{+}(s)$.
        \begin{definition}
          A path $p \in \gnpaths(s, L_n)$ is referred as a \emph{$\cal{G}$-augmenting} path for the vertex $s$
            if $\delta_{max}=\min\limits_{e \in p}(\flowpos(e) - \cal{G}(e)) > 0$.
        \end{definition}
        The latter means that the preflow $\cal{G}$ can be increased along this path by some constant $\delta \leq \delta_{max}$
          in such a fashion that the preflow $\cal{G}(e)$ does not exceed the flow $\flowpos(e)$ for each $e \in p$.
        Obviously, for a $\cal{G}$-augmenting path $p$ one has $\cal{G}(e) < \flowpos(e)$ for each $e \in p$.
        \begin{definition}
          The edge $e$ is called a \emph{$\cal{G}$-saturated} edge if $\cal{G}(e) = \flowpos(e)$.
          The path $p$ is called \emph{$\cal{G}$-saturated} if there exists such an edge $e \in p$ that $e$ is $\cal{G}$-saturated.
        \end{definition}
        For each path $p \in \gnpaths(s, L_n)$ we are going to raise the preflow $\cal{G}$ at every edge of the path $p$
        by some $\delta > 0$.
        We call the described procedure \emph{relaxation} of the vertex $s$.
        Strictly speaking, the {relaxation}
          of the vertex $s$ is a reduction of the excess of the flow $\cal{G}$ in the vertex $s$.
        Our purpose from now on is to relax the vertex $s$, transforming it into a $\cal{G}$-preserving vertex.

        On one hand, if $\delta_{max} \geq -d_{\cal{G}}(s)$, then \emph{pushing} the flow precisely $-d_{\cal{G}}(s)$ along the path $p$
          is enough to make the vertex $s$ preserving $\cal{G}$
          (\emph{pushing} the flow $C \in \mathbb{R}$ along the path $p$ means raising the flow $\cal{G}$ at each edge $e \in p$:
          $\cal{G}(e) \to \cal{G}(e) + C$).

        On the other hand, if $\delta_{max} < -d_{\cal{G}}(s)$, then pushing the flow $\delta_{max}$ along the path $p$ is not enough
          in order to make the vertex $\cal{G}$-preserving.
        In this case we push the flow $\delta_{max}$ along the path $p$ and look for the other augmenting paths in $\gnpaths(s, L_n)$.
        We will repeat the outlined procedure until we finally relax the vertex $s$.
        \begin{prop}
          Any vertex $s$ will be $\cal{G}$-relaxed in a finite number of steps.
        \end{prop}
        \begin{proof}
          Firstly, we have a finite graph $G_{n}$, and so there are only a finite number of paths in the $\gnpaths(s, L_n)$.
          Secondly, each time we push the flow along the augmenting path $p$ either we successfully relax the vertex $s$
            or we saturate at least one edge lying on the chosen path $p$.
          Since the number of the edges is finite, the procedure halts after a finite number of steps.

          Therefore, it happened that there are no augmenting paths in $\gnpaths(s, L_n)$.
          Examine the set of vertices $S \subseteq G_n$ that are reachable by non-saturated paths from the vertex $s$.
          We have $S \cap L_n = \emptyset$, hence $S \subseteq G_{n-1}$.
          Now consider the following subgraph $U = \bigcup\big\{p \in \gnpaths(s, S) \mid \text{$p$ is not $\cal{G}$-saturated}\big\}$.
          Similarly, the subgraph $U$ does not contain any vertices from $L_n$.
          \begin{lemma}
            Any edge $e \in \eout_G(U)$ is saturated.
          \end{lemma}
          \begin{proof}
            Due to the absence of the positive flow cycles, we have $\root \not \in U$.
            Since for each non-root vertex in $G_{k-1}$ we included every outgoing edge into $G_k$,
              the edge $e$ belongs to $G_n$.
            If $e$ is not saturated, then its ending vertex $\ter(e)$ lies in $S$.
            Therefore, the edge $e$ belongs to $U$.
            We arrived to contradiction.
          \end{proof}
          For a vertex $u \in S$ we see that $d_{\cal{G}}(u) = d^{+}_{\cal{G}}(u) - d^{-}_{\cal{G}}(u) \leq 0$, because
            the only $\cal{G}$-deficient vertex in the graph $G_n$ is the root vertex $\root$, which is not in $U$.
          Also for the vertex $s$ itself we have a strict inequality $d_{\cal{G}}(s) < 0$ since it is not relaxed by our assumption.
          It is clear now that the quantity $D = \sum_{u \in U} d_{\cal{G}}(u)$ is strictly less than zero.
          On the other hand,
          \[
            \begin{split}
              D &= \mathlarger\sum_{u \in U} \Big(d^{+}_{\cal{G}}(u) - d^{-}_{\cal{G}}(u)\Big)
              = \mathlarger{\smashoperator{\sum_{\eout_{G_n}(U)}}} \cal{G}(e) - \mathlarger{\smashoperator{\sum_{\ein_{G_n}(U)}}} \cal{G}(e) = \\
              &= \mathlarger{\smashoperator{\sum_{\eout_{G}(U)}}} \flowpos(e) - \mathlarger{\smashoperator{\sum_{\ein_{G_n}(U)}}} \cal{G}(e)
                \overset{\ref{p2}}{\geq} \mathlarger{\smashoperator{\sum_{\eout_G(U)}}} \flowpos(e) - \mathlarger{\smashoperator{\sum_{\ein_{G_n}(U)}}} \flowpos(e) 
                \overset{G_n \subset G} \geq \\
              &\geq \mathlarger{\smashoperator{\sum_{\eout_G(U)}}} \flowpos(e) - \mathlarger{\smashoperator{\sum_{\ein_{G}(U)}}} \flowpos(e) = 0.
            \end{split}
          \]
          Thus we deduce that $0 > D \geq 0$, and that is a contradiction.
        \end{proof}
        It follows that we can relax all the vertices in $\eback_n$.
        Conclusively, we set the preflow $\flowposn{n}$ equal to the $\cal{G}$ that we got after the relaxations.
        \begin{prop}
          The properties~\ref{p1}--\ref{p5} are satisfied for the resulting preflow $\flowposn{n}$ on the network $\net_n = (G_n, \len, \root)$.
        \end{prop}
        \begin{proof}
          The monotonicity property~\ref{p1} is true since during the relaxation step we
            only increased the flow on the edges of the graph $G_n$.

          The boundedness property~\ref{p2} is also true since we chose the flow values in such a way
            that the bound $\flowposn{n}(e) \leq \flowpos(e)$ holds for each edge $e\in G_n$.

          The root vertex was not affected during the step $n > 1$, so~\ref{p3} holds.

          The flow $\cal{G}$ was redesigned in such a fashion that each non-root vertex $v \in G_{n} \setminus 
            \left(L_n \cup \{\root\}\right)$ became $\cal{G}$-preserving.
          This checks the property~\ref{p4}.

          Finally, due to the property~\ref{posprop} and the fact that
          any vertex in $L_n$ has only incoming edges in the graph $G_n$, the property~\ref{p5} holds.
        \end{proof}
        The flow construction is completed.

        \medskip
        \subsection{Estimates on the constructed preflow}
          For abbreviation, let $\phi_n$ stand for $\phi_n(v) = \phi_{G_n}(v)$.
          Likewise, we define the functions $d_n$, $d^{-}_n$, $d^{+}_n$.
          \begin{lemma}
            For each $G_n$, $n > 1$,
            \begin{equation}
              \label{main-ineq}
              \sum_{v\in L_n}d_{n}^{-}(v) \phi_n(v) \leq \sum_{e \in E_n} \len(e) \flowpos(e).
            \end{equation}
          \end{lemma}
          \begin{proof}
            Set
            \[
              LS = \sum_{v \in L_n} d_n^{-}(v) \phi_n(v) = \sum_{v \in V_n} d_n^{-}(v) \phi_n(v).
            \]
            The last equality holds since $\phi_n(\root) = 0$ and $d_n(u) = 0$
              for any $u \in G_{n} \setminus \left(L_n \cup \{\root\}\right)$,

            Next we regroup the summation in order to have a sum over the edges of $G_n$, not vertices.
            Recall that $d_n^{+}(v) = 0$ for each vertex $v \in L_n$,
              hence $d_n(v) = -d_n^{-}(v)$ whenever $v$ is in the $n$-th layer.
            \begin{align*}
              LS = \sum_{v \in V_n} d_n^{-}(v) \phi_n(v) = -\sum\limits_{v \in V_n} d_n(v) \phi_n(v) =\\
              =-\sum_{v \in V_n} \left(\sum_{\eout(v)} \flowposn{n}(e) - \sum_{\ein(v)}\flowposn{n}(e) \right) \phi_n(v)=\\
              =\sum_{(vu) \in E_n} \flowposn{n}(vu) \big(\phi_n(u) - \phi_n(v)\big).
            \end{align*}
            For any $(vu) \in E_n$ we have $\phi_n(u) \leq \phi_n(v) + \len(vu).$ Therefore,
            \[
              LS\leq\sum_{(vu) \in E_n} \flowposn{n}(vu) \len(vu)= \sum_{e \in E_n} \flowposn{n}(e) \len(e) \overset{\ref{p2}}{\leq} \sum_{e \in E_n} \flowpos(e) \len(e).
            \]
            Inequality~\eqref{main-ineq} is proved.
          \end{proof}

          In what follows we prove that there is a vertex in any layer such that its root-distance is bounded by the mass of the flow.
          \begin{prop}
            \label{prop36}
            For each $n > 0$ there exists a vertex $r_n \in L_n$ such that $\phi_n(r_n) \leq \lvert\flowpos\rvert$.
          \end{prop}
          \begin{proof}
            We continue to write $LS$ for the left hand side of the inequality~\eqref{main-ineq}.
            Then
            \[
              LS = \sum_{v \in L_n}d^{-}_n(v) \phi_n(v) \geq \min_{v \in L_n}\phi_n(v) \sum_{v \in L_n} d^{-}_n(v).
            \]
            Since $\sum_{v \in L_n} d^{-}_n(v)$ is the total amount of the flow incoming to the deepest layer $L_n$, and the only
              active vertex in $G_n$ is the root, which produces the flow whose size is equal to $1$, we have $\sum_{v \in L_n} d^{-}_n(v) = d^{+}_n(\root) = 1.$
            Consequently, $LS \geq \min_{L_n} \phi_n(v).$

            Furthermore, combining the latter with~\eqref{main-ineq}, one gets
              \[
              \min_{L_n} \phi_n(v) \leq \sum_{e \in E_n} \len(e) \flowpos(e) \leq \lvert\flowpos\rvert.
              \]
            Since $L_n$ contains a finite number of vertices, there exists a vertex $r_n \in L_n$ such that
              $\phi_n(r_n) \leq \lvert\flowpos\rvert.$
          \end{proof}
          \begin{prop}
            There exists a ray $r=(r_1 r_2\dots)$ such that $\sum\limits_{k=1}^\infty \len(r_kr_{k+1}) < \infty$.
          \end{prop}
          \begin{proof}
            We will seek for the ray $r$ in the graph $\gstar = \bigcup G_n.$
            We abbreviate $\phi_{\gstar}(v)$ to $\phistar$.
            Consider the paths $p$, starting on the root, such that $\phistar(p) \leq \lvert\flowpos\rvert$.
            Denote this set of paths by $\cal{P}^*$.

            For each $k > 0$ we are going to present a vertex $r_k$ such that there are an infinite number of paths
              in $\cal{P}^*$ starting on $(r_1 r_2 \dots r_k)$.

            At the first step we set $r_1 = \root$.
            By proposition~\ref{prop36}, for each layer $L_n$ there is a vertex $v_n \in L_n$ such that
              $\phi_n(v_n) \leq \lvert\flowpos\rvert$.
            Since $\phistar(v_n) \leq \phi_n(v_n)$, and there are an infinite number of layers,
              $\cal{P}^*$ contains an infinite number of paths.

            Now assume we have already constructed the path $p_k = (r_1 r_2\dots r_k)$ such that
              there are an infinite number of paths in $\cal{P}^*$ starting with the path $p_k$.
            The graph $\gstar$ is locally finite and in particular the degree of $r_k$ is finite.
            Since the number of the paths from $\cal{P}^*$ starting with $p_k$ is infinite by the assumption,
              we are able to choose $r_{k+1} \in \eout_{G^*}(r_k)$ such that the number of paths in $\cal{P}^*$ starting from
              the path $p_{k+1}= (p_k r_k)$ is infinite as well.

            We constructed the ray $r = (r_1 r_2 \dots) \subseteq \gstar$.
            The total length of the ray $\sum \len(r_k r_{k+1})$ is finite
              since the partial sums are bounded by $\lvert \flowpos \rvert$.
          \end{proof}
          We found a finite-length ray in the network $\net$.
          Hence, the implication is proved as well as the theorem itself.
      \end{proof}

      \medskip
      \subsection{Main theorem}
        Now we present the proof of Theorem~\ref{thm-main} using Theorems~\ref{thm-graph-eq} and~\ref{ray-thm}.
        \begin{proof}[Proof of Theorem~\ref{thm-main}]
          The system $\fsys$ belongs to B-class, so we can apply Theorem~\ref{thm-graph-eq}, which asserts
            that $\fsys$ admits the rank one density property if and only if
            the corresponding network $\net(\fsys)$ is $\flow$-preserving for any
            finite-mass flow $\flow$.
          We rephrase the last conclusion in terms of \emph{oriented} networks.
          The system $\fsys$ admits the rank one density property if and only if
            the oriented network $\net(\fsys)$ is $\flowpos$-preserving for any
            finite-mass flow $\flowpos$.
          After that we utilize Theorem~\ref{ray-thm}, which states that there exists such
            flow $\flowpos: E \to \mathbb{R}^{+}$ with a finite mass that is not preserved in
            the network $\net(\fsys)$, if and only if there exists a ray $r^{*} \subseteq \net$, originating in the vertex $\source$,
            whose total length is finite $\sum_{k=1}^\infty \len(r^{*}_k r^{*}_{k+1}) < \infty$.
          Combining these two statements we get that $\fsys$ is \emph{not} rank one dense if and only if
            there is a finite-length ray $r^*$ in the network $\net(\fsys)$.
          It is obvious that $r^*$ cannot include infinite number of the edges incident to the vertices
            $\source$ and $\sink$, since such edges have the length exactly one.
          Therefore, such $r^*$ exists in $\net(\fsys)$ if and only if there is a ray $r=\seq{r}$ in the original bipartite graph $B(\fsys)$ such that
            the series $\sum_{k=1}^\infty\lvert\wt(r_k, r_{k+1})\rvert$ converges.
          The theorem is proved.
        \end{proof}

  \bigskip
  \section{Acknowledgements}
    The author wishes to express his gratitude to Dmitry Yakubovich for suggesting the construction of $M$-bases using the infinite bipartite graphs.
    The author gratefully acknowledges the many helpful suggestions of Anton Baranov during the preparation of the paper.

\end{document}